\documentclass[leqno,11pt]{amsart}
\usepackage{amssymb,amsmath,multirow,graphics,multicol}

\newtheorem{theorem}{Theorem}[section]

\newtheorem{proposition}{Proposition}[section]
\newtheorem{corollary}{Corollary}[section]

\theoremstyle{remark}
\newtheorem{remark}{Remark}[section]
\theoremstyle{remark}

\theoremstyle{definition}
\newtheorem{definition}{Definition}[section]

\numberwithin{equation}{section}

\def\<{\left < }
\def\>{\right >}
\def\({\left ( }
\def\){\right )}

\begin{document}

\thispagestyle{plain}

\title[Sharp Growth estimates for warping functions]{Sharp growth estimates for warping functions in multiply warped product manifolds}
\author[B.-Y. Chen  and S. W. Wei]{Bang-Yen Chen  and Shihshu Walter Wei*}
\date{}
\thanks{*The second author  is supported in part by NSF (DMS-1447008).}

\maketitle

\begin{abstract} By applying an average method in PDE, 
we obtain a dichotomy between ``constancy'' and ``infinity" of the warping functions on complete noncompact Riemannian manifolds for an appropriate isometric immersion of a multiply
warped product manifold $N_1\times_{f_2} N_2 \times \cdots \times _{f_k} N_k\, $ into a Riemannian manifold. 

Generalizing the earlier work of the authors in [9], we establish sharp inequalities between
 the mean curvature of the immersion and the sectional curvatures  of  the  ambient  manifold  under the influence of quantities  of  a  purely  analytic  nature  (the growth of the warping functions). Several applications of our growth estimates are also presented. \\[0.2cm]
\end{abstract}

\label{first}

\section{\bf Introduction}

Warped products play very important roles in Differential Geometry and Physics. Examples of warped product include Riemannian manifolds of constant curvature and the best relativistic model of the Schwarzschild space-time that describes the out space around a mass star or a black hole. 

In \cite {CW}, B.-Y. Chen and S. W. Wei obtained the following
necessary condition for an arbitrary  isometric immersion of a
warped product $N_1 \times _f N_2$ into a Riemannian $m$-manifold $\tilde M^m_c$ with sectional curvatures bounded from above by a constant $c$, generalizing the work of B.-Y.  Chen in \cite {C} on warped product submanifolds in a Riemannian manifold $R^m(c)$ of constant sectional curvature $c$.
\medskip

\noindent
{\bf Theorem A.} \cite [Theorem 3.1] {CW} {\it For any isometric immersion $\phi : N_1 \times _f N_2 \to \tilde M^m_c$   from a warped product $N_1 \times _f N_2$ into  a Riemannian $m$-manifold $\tilde M^m_c$ with sectional curvatures bounded from above by a constant $c$, the warping function $f$ satisfies
\begin{align}\label{1.1} -\frac{(n_1 + n_2)^2}{4 n_2} H^2 - n_1c \le  \frac {\Delta f}{f},\end{align}
where $n_1= \dim N_1$ and $n_2=\dim N_2$,    $H^2 = \langle H, H \rangle$ is the squared mean curvature of $\phi\, ,$ and $\Delta f$ is the  Laplacian of $f$   on $N_1$ $($defined as the divergence of the gradient vector field of $f$, cf. \eqref{2.4}$)$.   

The equality sign \eqref{1.1} holds if and only if  $\phi$ is a mixed totally geodesic immersion with $\operatorname{trace} h_1 =  \operatorname{trace} h_2$, where $h_1$ and $h_2$ are the restriction of the second fundamental form $h$ of $\phi$ restricted to $N_1$ and $N_2$, respectively, and
at each point $p=(p_1,p_2) \in N , $ $c$ satisfies $c = K(u, v) = \max K(p)$, for every unit vector $u \in T^1 _{p_1} N_1$ and every unit vector $v \in T^1 _{p_2} (N_2).$  }
\vskip.1in

On the other hand, the second author extended  in \cite {W} the scope of $L^q$ or $q$-integrable functions on complete noncompact Riemannian manifolds  to functions with ``\emph{$p$-balanced}" growth depending on $q$, and introduced the concepts of their counter-part to  ``\emph{$p$-imbalanced}" growth (cf. Definition \ref{def2.1}). By coupling these growth estimates with the above inequality \eqref{1.1}, Chen and Wei establish in \cite{CW} some sharp inequalities between
quantities of a geometric nature (the mean curvature of the immersion, the sectional
curvatures  of  the  ambient  manifold)  and  quantities  of  a  purely  analytic  nature  (the
growth of the warping function).

\vskip.1in

\noindent
{\bf Theorem B.} \cite {CW} {\it  If $f$ is nonconstant and \emph{$2$-balanced} for some $q > 1$,
then for every Riemannian $n_2$-manifold $N_2$ and every isometric
immersion $\phi$ of the warped product $N_1\times_f N_2$ into any
Riemannian manifold $\tilde M^m_c$ with $c\leq 0$, the mean
curvature $H$ of $\phi$ satisfies
\begin{align}\label{1.2}H^2 > \frac{4n_1 n_2|c|}{(n_1 +
n_2)^2}\end{align} at some points.}
\vskip.1in

 Hence we immediately find a dichotomy between
``constancy'' and ``infinity" ($2$-imbalanced) of the warping functions on complete
noncompact Riemannian manifolds for an appropriate isometric
immersion: \vskip.1in

\noindent
{\bf Corollary A.} \cite {CW}  {\it  Suppose the squared mean curvature of the
isometric immersion $\phi:N_1\times_f N_2\to \tilde M^m_c$
satisfies
\begin{align}\label{1.3}H^2 \leq \frac{4n_1 n_2|c|}{(n_1 + n_2)^2}\end{align} everywhere on $N_1\times_f
N_2$. Then the warping function $f$ is either a constant or it has \emph{$2$-imbalanced} growth for
every $q > 1$.}
\vskip.1in

Applications of these new inequalities are also presented, among which there are some results on the nonexistence of isometric minimal
immersions between certain types of Riemannian manifolds:
\vskip.1in

\noindent
{\bf Theorem C.} \cite {CW}  {\it  
Suppose $q > 1$ and the warping function $f$ is \emph{$2$-balanced}.
 If $N_2$ is compact,  then there does not exists an isometric minimal immersion from  $N_1 \times _f N_2 $ into any Euclidean space.}
\vskip.1in

A Riemannian manifold is said to be {\it  negatively curved} (resp., {\it non-positively curved}\/) if it has negative (resp., {\it non-positively curved}\/)  sectional curvatures.
\vskip.1in

\noindent
{\bf Corollary B.} \cite {CW}  {\it  
If $f$ is an $L^q$ function on $N_1$ for some $q > 1$, then for any Riemannian manifold $N_2$ the  warped product  $N_1\times_f N_2$ does not admit any isometric minimal immersion into any non-positively curved Riemannian manifold.}
\vskip.1in

For further extension, let $N = N_1 \times \cdots \times N_k$ denote the Cartesian product of $k$ Riemannian manifolds $(N_1, g_1) \cdots,  (N_k, g_k)\, ,$ and $\pi _i : N \to N_i$ be the canonical projection of $N$ onto $N_i\, , 1 \le i \le k.$ If $f_2, \cdots f_k: N_1 \to \mathbb R^{+}$ are smooth positive-valued functions, then   
\[ g = \pi_1^{\ast} g_1 + \sum _{i=2}^{k} (f_i \circ \pi_1)^2\pi_i^* g_i\] 
defines a Riemannian metric on $N\, ,$ called  multiply warped product metric. The product manifold $N$ endowed with $g$ is denoted by $ N = N_1\times_{f_2} N_2 \times \cdots \times _{f_k} N_k\, .$ 

Denote by $\operatorname{trace}  h_i$ the trace of the second fundamental form $h$ of  $N = N_1 \times \cdots \times N_k$ into a Riemannian manifold restricted to $N_i\, .$ 

B.-Y. Chen and F. Dillen proved in \cite{CD} the following.
\vskip.1in

\noindent
{\bf Theorem D.}  {\it  Let $\phi : N_1\times_{f_2} N_2 \times \cdots \times _{f_k} N_k \to M$ be an isometric 
immersion  of a multiply warped product $N = N_1\times_{f_2} N_2 \times \cdots \times _{f_k} N_k $ into an arbitrary Riemannian manifold $M\,$. Then we have
\begin{equation}\label{1.4} - \frac{n^2(k-1)}{2k}H^2- n_1(n - n_1) \max \tilde{K} \leq \sum _{j=2}^{k} n_j  \frac{\Delta f_j}{f_j} ,\end{equation}
where $n=\sum_{i=1}^k n_i$ and $\max\tilde{K}(p)$ denotes the maximum of the sectional curvature function of the ambient space $M$ restricted to $2$-plane sections of the tangent space $T_pN$ of $N$ at $p = (p_1,\dots, p_k)$.

The equality sign of \eqref{1.4} holds identically if and only if the following two conditions hold 

\begin{enumerate}
\item[{\rm i)}] $\phi$ is mixed totally geodesic such that $\;\operatorname{trace}\, h _1 = \cdots = \operatorname{trace}\, h _k\, ;$ 

\item[{\rm ii)}] At each point $p \in N$, the sectional curvature function $\tilde K$ satisfies $\tilde K(u, v) = \max \tilde K(p)$ for every $u \in T^1 _{p_1} N_1$ and every  $v \in T^1 _{p_2, \dots p_k} (N_2 \times \cdots N_k)\,$.
\end{enumerate}}
\vskip.1in

One main purpose of this article is to prove the following theorem which extends Theorem B, in particular, inequality \eqref{1.2} to arbitrary isometric immersions of  multiply warped product manifolds into an arbitrary Riemannian manifold. 
 
\begin{theorem} \label{T:1.1}
If for each $j,\, 2\leq j\leq k,$ $f_j$ is nonconstant and \emph{$2$-balanced} with some $q_j > 1$,
then, for any multiply warped product $N = N_1\times_{f_2} N_2 \times \cdots \times _{f_k} N_k$ in a
Riemannian manifold $M$, the mean curvature $H$ of $N$ in $M$ satisfies
\begin{equation}\label{1.5}H^2 > \frac{- 2kn_1(n - n_1) }{n^2(k - 1)}\max \tilde{K}\end{equation} at some points, where $\max \tilde{K}$ is defined in Theorem D. 

In particular, if each $f_j $ is nonconstant and in $L^{q_j}$ for some $q_j > 1$, then \eqref{1.5} holds at some points.\end{theorem}

In particular, if $M$ is a Riemannian manifold of constant sectional curvature $c\leq 0$, then Theorem \ref{T:1.1} reduces to the following.

\begin{theorem} \label{T:1.2}
If for each $j,\, 2\leq j\leq k,$ $f_j$ is nonconstant and \emph{$2$-balanced} with some $q_j > 1$,
then, for any multiply warped product $N = N_1\times_{f_2}\hskip-.02in  N_2 \times \cdots \times _{f_k}\hskip-.02in N_k$ in a Riemannian manifold $R^m(c)$ of constant sectional curvature $c \leq 0$, the mean curvature $H$ of $N$ in $R^m(c)$ satisfies
\begin{equation}\label{1.6}H^2 > \frac{2kn_1(n - n_1) }{n^2(k - 1)}c\end{equation} 
at some points.  

In particular, if each $f_j $ is nonconstant and in $L^{q_j}$ for some $q_j > 1$, then \eqref{1.6} holds at some points.\end{theorem}

Theorems \ref{T:1.1} and \ref{T:1.2} are sharp and inequalities \eqref{1.5} and \eqref{1.6} are optimal. For details, we refer to Remark \ref{R:3.1}, Example 3.1, Example 3.2, and Remark \ref{R:3.2}.

In views of Theorem \ref{T:1.1}, we give the following dichotomy.

\begin{theorem} \label{T:1.3} Suppose the squared mean curvature of the isometric immersion of a multiply warped product $N = N_1\times_{f_2} N_2 \times \cdots \times _{f_k} N_k$ into a Riemannian manifold satisfies
\begin{align}\label{1.7}H^2 \leq \frac{- 2kn_1(n - n_1) }{n^2(k - 1)}\max \tilde{K}\end{align} everywhere on $N$. Then there exists an integer $i, 2 \le i \le k$ such that either  the warping function $f_i$ is a constant or  $f_i$ has \emph{$2$-imbalanced} growth for every $q_i > 1$. 
\end{theorem}

Some other applications of Theorem \ref{T:1.1} are the following. 

\begin{corollary} \label{C:1.1} If for each $j, 2 \le j \le k$, $f_j$ is nonconstant and \emph{$2$-balanced} for some $q_j > 1$,
then there does not exist a minimal immersion of any multiply warped product $N = N_1\times_{f_2} N_2 \times \cdots \times _{f_k} N_k$ into a
Riemannian manifold whose maximum sectional curvature is nonpositive. 

In particular,  if each $f_j $ is nonconstant and in $L^{q_j}$ for some $q_j > 1$, then there does not exist a minimal immersion of any multiply warped product $N = N_1\times_{f_2} N_2 \times \cdots \times _{f_k} N_k$ into a Euclidean space.
\end{corollary}

Applying the growth estimates in Theorem \ref{T:2.1} and the average method in PDE in Proposition \ref{P:2.1}, we have the following Liouville property and characterization results.

\begin{corollary} \label{C:1.2}
Suppose the squared mean curvature of the isometric immersion $\phi$ of a multiply warped product $N = N_1\times_{f_2} N_2 \times \cdots \times _{f_k} N_k$ into a complete, simply-connected Riemannian manifold $R^m(c)$ of constant sectional curvature $c$
satisfies \eqref{1.7} everywhere on $N$.  If for each $j, 2 \le j \le k$, $f_j$ is \emph{$2$-balanced} for some $q_j > 1$, then we have:
\begin{enumerate}
\item[{\rm (1)}] Every  warping function $f_j, 2 \le j \le k$ is constant. 

\item[{\rm (2)}] The isometric immersion $\phi$ is a minimal immersion into a Euclidean space. 

\item[{\rm (3)}] The isometric immersion $\phi$ is a warped product immersion. 
\end{enumerate}\end{corollary}

\begin{corollary} \label{C:1.3} Let each $f_j\, ,$ $2 \le j \le k$  be  \emph{$2$-balanced} for some $q_j > 1\, .$ Then we have:
\begin{enumerate}
\item[{\rm (1)}] Every multiply warped product $N = N_1\times_{f_2} N_2 \times \cdots \times _{f_k} N_k$ does not admit  an isometrically minimal immersion into any
Riemannian manifold of negative sectional curvature.
\item[{\rm (2)}] If $N_k$ is compact,  then  $N = N_1\times_{f_2} N_2 \times \cdots \times _{f_k} N_k$ does not admit  an isometrically minimal immersion into a Euclidean space. 
\end{enumerate} \end{corollary}

We state a special case of Corollary \ref{C:1.3} as the following.

\begin{corollary} \label{C:1.4}   If each $f_j\, ,$ $2 \le j \le k,$ is  in $L^{q_j}$ for some $q_j > 1$, then we have:

\begin{enumerate}
\item[{\rm (1)}] Every multiply warped product $N = N_1\times_{f_2} N_2 \times \cdots \times _{f_k} N_k$ does not 	admit  an isometrically minimal immersion into any negatively curved Riemannian manifold.

\item[{\rm (2)}]  If $N_k$ is compact,  then $N = N_1\times_{f_2} N_2 \times \cdots \times _{f_k} N_k$ does not admit  an isometrically minimal immersion into a Euclidean space. 
\end{enumerate}\end{corollary}

A map $$\psi : N_1\times_{f_2} N_2 \times \cdots \times _{f_k} N_k \to M_1\times_{{\rho}_2} M_2 \times \cdots \times _{{\rho}_k} M_k$$ between two multiply warped product manifolds $N_1\times_{f_2} N_2 \times \cdots \times _{f_k} N_k$ and $M_1\times_{{\psi}_2} M_2 \times \cdots \times _{{\psi}_k} M_k$ is said to be a {\it warped product immersion} if $\psi$ is given by $\psi (x_1, \cdots , x_k) = (\psi _1(x_1), \cdots , \psi _k(x_k))$ is an isometric immersion, where $\psi_i : N_i \to M_i, i=2, \cdots , k$ are isometric immersions, and $f_i =  \rho_i \circ \psi_1 : N_1 \to {\mathbb R}^{+}$ for $i=2, \cdots , k\, .$

By applying Theorem D, Proposition \ref{P:2.1} and Theorem E (N\"olker's Theorem), we have

\begin{corollary} \label{C:1.5} If for each $j, 2 \le j \le k$, $f_j$ is \emph{$2$-balanced} for some $q_j > 1$,
then  every  isometric minimal immersion of a multiply warped product $N = N_1\times_{f_2} N_2 \times \cdots \times _{f_k} N_k$ into a Euclidean space  is a warped product immersion. 
\end{corollary}

The technique used in this article is to apply the Average Method in PDE in Proposition \ref{P:2.1} and the Growth Estimates in Theorem \ref{T:2.1} to study multiply warped products. In contrast to {\it an extrinsic average variational method in the calculus of variations} \cite {W2, W3,  HW}, where the sum of analytic quantities is strictly negative, {\it the average method in PDE} in this article deals with the nonnegative sum of analytic quantities (cf. Remark 2.1). 

The techniques used in this article are sufficient general to apply to multiply warped product  manifolds totally real isometrically immersed into complex space forms, as well as into quaternionic space forms. We also use the same technique for multiply warped product manifolds to treat doubly warped product manifolds in the last section.

\section{Preliminaries.}

Let $N$ be a Riemannian $n$-manifold isometrically immersed in a Riemannian $m$-manifold $\tilde M^m$.   We choose a local field of orthonormal frame
$e_1,\ldots,e_n,e_{n+1},\ldots,e_m$ in $\tilde M^m$ such that, restricted
to  $N$, the vectors $e_1,\ldots,e_n$ are tangent to $N$
and $e_{n+1},\ldots,e_m$ are normal to $N$.

For a submanifold $N$ in $\tilde M^m$, let $\nabla$ and ${\tilde\nabla}$ denote the Levi-Civita connections of $N$ and $\tilde M^m$, respectively. The Gauss and Weingarten formulas are then given respectively by (see, for instance, \cite{c1,C3})
\begin{align}\label{2.1} &{\tilde \nabla}_{X}Y=\nabla_{X} Y +h(X,Y),\\ & \label{2.2}{\tilde\nabla}_{X}\xi =
-A_{\xi}X+D_{X}\xi \end{align} for  vector
fields $X,Y$ tangent to $N$ and  $\xi$ normal to $N$, where $h$ is the second fundamental form, $D$ the normal connection, and $A$ the shape operator of the submanifold.
Let $\{h^r_{ij}\}$, $i,j=1,\ldots,n;\,r=n+1,\ldots,m$, denote the coefficients of the second fundamental form $h$ with respect to $e_1,\ldots,e_n,e_{n+1},\ldots,e_m$.

The mean curvature vector $\overrightarrow{H}$ is defined by
\begin{align}\label{2.3}\overrightarrow{H} = {1\over n}\,\hbox{\rm trace}\,h = {1\over n}\sum_{i=1}^{n}
h(e_{i},e_{i}), \end{align}
where $\{e_{1},\ldots,e_{n}\}$ is a local  orthonormal
frame of the tangent bundle $TN$ of $N$. The
 squared mean curvature is given by
$$H^2=\left<\right.\hskip-.02in\overrightarrow{H},\overrightarrow{H}\hskip-.02in
\left.\right>,$$
where $\<\;\,,\;\>$ denotes the inner product. A submanifold $N$  is called {\it minimal}  in $\tilde M^m$ if its mean curvature vector vanishes identically.

Let $P$ be a Riemannian $k$-manifold and $\{e_1,\ldots,e_k\}$ be an
orthonormal frame field on $P$. For a differentiable function
$\varphi$ on $P$, the {\it Laplacian} of $\varphi$ is defined  by the
divergence of the gradient of $\varphi$, or the trace of the
Hessian $\varphi,$ i.e.
\begin{align} \label{2.4} &\Delta\varphi=\sum_{j=1}^k \{e_je_j\varphi - (\nabla_{e_j}e_j)\varphi \}.\end{align}

A function $\varphi$ on $P$ is said to be {\it harmonic} (resp. {\it subharmonic} or {\it superharmo\-nic}\/) if we have $\Delta\varphi =0$ (resp. $\Delta\varphi\ge 0$ or  $\Delta\varphi\le 0$)  on $P$.

An isometric immersion $$\phi: N_1\times_{f_2} N_2 \times \cdots \times _{f_k} N_k \to M$$ of a multiply warped product $ N_1\times_{f_2} N_2 \times \cdots \times _{f_k} N_k $ into a Riemannian $m$-manifold $M$ is called {\it mixed totally geodesic\/} if its second fundamental form $h$ satisfies
$h(\mathcal D_i,\mathcal D_j)=0$ for any distinct $i,j \in \{1, \cdots, k\},$ where $\mathcal D_i $ denotes the distribution obtained from the vectors tangent to the horizontal lifts of $N_i\, .$  

We recall the following results for later use.

\vskip.1in
\noindent
{\bf Theorem E.} \cite[N\"olker's Theorem] {N} {\it Let $\phi: N_1\times_{f_2} N_2 \times \cdots \times _{f_k} N_k \to R^m(c)$ be an  isometric immersion into a Riemannian manifold $R^m(c)$ of constant sectional curvature $c\, .$ If $\phi$ is mixed totally geodesic, then locally $\phi$ is a warped product immersion} \smallskip
\vskip.1in

In the following, let us assume that $N_1$ is a noncompact complete Riemannian manifold  and  $B(x_0;r)$ denotes  the geodesic ball of radius $r$
centered at $x_0 \in N_1$.

We recall some notions from \cite{W}.
\vskip.1in

\begin{definition}\label{def2.1}
A  function  on $N_1$ is said to have {\rm $p$-balanced growth} $($or, simply,  {\rm is $p$-balanced}$)$ if it is one of the following:  $p$-finite,  $p$-mild,   $p$-obtuse,  $p$-moderate, and $p$-small; it has {\rm $p$-imbalanced growth}, or simply is {\rm $p$-imbalanced} otherwise. \end{definition}

\noindent Notice that the definitions of   ``\emph{$p$-{finite}}, \emph{$p$-mild},  \emph{$p$-obtuse},  \emph{$p$-moderate}, \emph{$p$-small}" and their counter-parts ``\emph{$p$-infinite}, \emph{$p$-severe},  \emph{$p$-acute}, \emph{$p$-immoderate}, \emph{$p$-large}" growth depend on $q$, and
$q$ will be specified in the context in which the definition is used. 

We have discussed their definitions in \cite [Definition 4.1-4.5] {CW}. For completeness we include them as follows (please see also \cite {W}).
\smallskip

\begin{definition} \label{def2.2}
A  function $f$  on $N_1$ is said to have {\rm $p$-{finite growth}} $($or, simply, {\rm is $p$-{finite}}$)$ if there exists $x_0 \in N_1$ such that
\begin{equation} \lim_{r \to \infty} \inf\frac{1}{r^p}\int_{B(x_0;r)} |f|^{q}dv < \infty ;   \label{7.1} \end{equation}
it has \emph{$p$-{infinite growth}} $($or, simply, \emph{is
$p$-infinite}$)$ otherwise.\end{definition}

\begin{definition}\label{def2.3}
A function $f$ has \rm{$p$-mild growth} $($or, simply, \rm{is $p$-mild}$)$ if there exists  $ x_0 \in N_1\, ,$ and a strictly increasing sequence of $\{r_j\}^\infty_0$ going to infinity, such
that for every $l_0>0$, we have
\begin{equation}\sum\limits_{j=\ell_0}^{\infty}\bigg(\frac{(r_{j+1}-r_j)^p}{\int_{B(x_0;r_{j+1})\backslash B(x_0;r_{j})}|f|^qdv}\bigg)^{\frac1{p-1}}=\infty \,;   \label{2.7}\end{equation}
and has \rm{$p$-severe growth} $($or, simply, \rm{is
$p$-severe}$)$ otherwise.\end{definition}

\begin{definition}\label{def2.4}
A function $f$ has \emph{$p$-obtuse growth} $($or, simply, \emph{is $p$-obtuse}$)$ if there exists $x_0 \in N_1$ such that for every $a>0$, we have
\begin{equation}\int^\infty_a\bigg( \frac{1}{\int_{\partial B(x_0;r)}|f|^qdv}\bigg)^\frac{1}{p-1}dr = \infty \, ;   \label{2.8}\end{equation}
and has \emph{$p$-acute growth} $($or, simply, \emph{is
$p$-acute}$)$ otherwise.
\end{definition}

\begin{definition}\label{def2.5}
A function $f$ has \emph{$p$-moderate growth} $($or, simply, \emph{is $p$-moderate}$)$ if there exist  $ x_0 \in N_1$, and
$F(r)\in {\mathcal F}$,such that \begin{equation}
\lim \sup _{r \to \infty}\frac {1}{r^p F^{p-1} (r)}\int_{B(x_0;r)} |f|^{q}dv < \infty .   \label{2.9}
\end{equation}
And it has \emph{$p$-immoderate growth} $($or, simply, \emph{is
$p$-immoderate}$)$ otherwise,  where
\begin{equation} {\mathcal F} = \{F:[a,\infty)\longrightarrow
(0,\infty) |\int^{\infty}_{a}\text{\small$
\frac{dr}{rF(r)}$}=+\infty \ \ for \ \ some \ \ a \ge 0 \}\, .
\label{2.10}\end{equation} {\rm (Notice that the functions in
{$\mathcal F$} are not necessarily monotone.)} 
\end{definition}

\begin{definition}\label{def2.6}
A function $f$ has \emph{$p$-small growth} $($or, simply, \emph{is $p$-small}$)$ if there exists $ x_0 \in N_1\, ,$ such that for every $a>0\, ,$we have
\begin{equation}\int^\infty_a\bigg( \frac{r}{\int_{ B(x_0;r)}|f|^qdv}\bigg)^\frac{1}{p-1}dr = \infty ;   \label{2.11}
\end{equation}
and  has \emph{$p$-large growth} $($or, simply, \emph{is
$p$-large}$)$ otherwise.
\end{definition}

We recall  the following result from \cite {CW} for later use.

\begin{theorem} [Warping Function Growth Estimates]\label{T:2.1} Let $N_1$ be a noncompact complete  Riemannian manifold and $f: N_1 \to \mathbb R^{+}$ be a $C^2$ positive function satisfying ${\Delta f}/{f} \ge 0$ on $N_1\, .$ Then either $f$ is constant or $f$ is  $2$-imbalanced for every $q > 1$.
\end{theorem}
\vskip.05in

\noindent{\it Proof}. Follow exactly the proof of Theorems 4.1, 4.2, 4.3, 4.4 and 4.5 in \cite [p.586-590] {CW}
 and use Definition \ref{def2.1}, the assertion follows. \hfill $\square$

We also need the following result.

\begin{proposition}[An Average Method in PDE] \label{P:2.1} Let  $c_2, \cdots, c_k $ be $k-1$ positive constants and let $f_2, \cdots , f_k$ be positive-valued functions defined on a complete noncompact manifold $N_1$ such that $\sum _{j=2}^k c_j  {\Delta f_j}/{f_j} \ge 0$, Then we have:

\begin{enumerate}
\item[{\rm (1)}] There exists an integer $i,\, 2\le i \le k$, such that either $f_i$ is a constant or $f_i$ is  $2$-imbalanced for every $q_i > 1$.

\item[{\rm (2)}]  If each $f_j\, ,$ $2 \le j \le k$, is  $2$-balanced for some $q_j > 1$, then all of $f_2, \cdots , f_k$ are constant functions.
\end{enumerate}\end{proposition}
\vskip.1in

\noindent{\it Proof.}  If  $\sum _{j=2}^k c_j  {\Delta f_j}/{f_j} \ge 0$ holds, then there exists at least $i,\, 2 \le i\le k$, such that ${\Delta f_i}/{f_i} \ge 0$ holds. Or $\sum _{j=2}^k c_j {\Delta f_j}/{f_j} < 0\, ,$ contradicting to the hypothesis. Therefore statement (1) of this proposition follows from Theorem \ref{T:2.1}.

For statement (2), it follows from the assumptions that $f_i$ is $2$-balanced. Hence statement  $(1)$ implies that  $f_i$ is constant. So we have $\Delta f_i = 0$ and $$\sum _{j\ne i}  c_j   \frac{\Delta f_j}{f_j} = \sum _{j=2}^k c_j \frac {\Delta f_j}{f_j} \ge 0.$$
Now, by applying statement $(1)$ again to $\sum _{j\ne i}  c_j   {\Delta f_j}/{f_j}  \ge 0$, we can find the second constant warping function $f_{i^{\prime}}$ such that $ \sum _{j\ne i, i^{\prime}}  c_j   {\Delta f_j}/{f_j}\ge 0.$ Now, using the same method iteratively,  we conclude that all of $f_2, \cdots , f_k$ are all constant. This proves statement (2).
\hfill $\square$

\begin{remark} \label{R:2.1} {\rm The average method given in Proposition \ref{P:2.1} is in contrast to {\it an extrinsic average variational method in the calculus of variations} \cite {W2,W3,HW}, where the sum of analytic quantities, the second variation formulas of functionals such as the mass, $p$-energy, or Yang-Mills functional (over a set of distinguished variation vector fields) is {\it strictly negative}.  Our average method in PDE in Proposition \ref{P:2.1} deals with the {\it nonnegative} sum of analytic quantities, the Laplacian of warping functions.}
\end{remark}

\section{\bf Proof of Theorem \ref{T:1.1}, Theorem \ref{T:1.3} and Corollaries  \ref{C:1.1} - \ref{C:1.5}.}

The proofs of these results are based on Theorem D via the Average Method in PDE given in  Proposition \ref{P:2.1} and the Warping Function Growth Estimates given in Theorem \ref{T:2.1}.

\vskip.05in
\noindent{\it Proof of Theorem \ref{T:1.1}}. Suppose contrary to \eqref{1.5}, i.e., there were an isometric immersion $\phi$ whose mean curvature $H$ satisfying
\begin{equation}\label{3.1}H^2 \leq \frac{- 2kn_1(n - n_1)\max \tilde{K} }{n^2(k - 1)}\end{equation}  everywhere on $N\, .$
This would imply by multiplying both sides of \eqref{1.7} by a positive number $ {n^2(k - 1)}/{2k}\, ,$ or equivalently
$$0 \le - \frac {n^2 (k - 1)}{2k} H^2 - n_1 (n - n_1) \max \tilde{K}. 
$$

On the other hand, Theorem D would imply
\begin{equation}\notag - \frac{n^2(k-1)}{2k}H^2- n_1(n - n_1) \max \tilde{K} \leq \sum _{j=2}^{k} n_j  \frac{\Delta f_j}{f_j} .\end{equation}
After combining this with \eqref{3.1} or its equivalent inequality, we find
\begin{equation} \label{3.2} 
 \sum _{j=2}^{k} n_j  \frac{\Delta f_j}{f_j}  \ge  - \frac {n^2 (k - 1)}{2k} H^2 - n_1 (n - n_1) \max \tilde{K} \ge 0.
 \end{equation}
 
 Now, after applying the Average Method in PDE stated in Proposition \ref{P:2.1}(1),  we would conclude from \eqref{3.2}  that 
 some $f_i,\, 2\leq i\leq k,$ could be constant {\it or} $f_i$ would be $2$-imbalanced for every $q_i > 1$, contradicting the assumption that
$f_i$ is nonconstant {\it and}  $2$-balanced for some $q_i > 1\, .$ Indeed, ``$f_i$ would be constant" contradicts ``$f_i$ is nonconstant" and 
``$f_i$ would be $2$-imbalanced for every $q_i > 1$." contradicts ``$f_i$ is $2$-balanced for some $q_i > 1$.
\vskip.1in

To prove the last assertion, we observed that every $L^q$ function has
\emph{$2$-{finite}}, \emph{$2$-mild}, \emph{$2$-obtuse},
\emph{$2$-moderate}, \emph{$2$-small} growth for the same $q$ (cf.
\cite [Proposition 2.3] {WLW}). For example, if $f$ defined on $N_1$ is in $L^q$, then $f$ is $2$-finite with respect to the same $q$. Indeed, 
there exists $x_0 \in N_1$ such that \eqref{7.1}, where $p=2$ holds:
\[ \begin{aligned}
\lim_{r \to \infty} \inf\frac{1}{r^2}\int_{B(x_0;r)} |f|^{q}dv & \le \lim_{r \to \infty} \inf\frac{1}{r^2}\int_{N_1} |f|^{q}dv\\
                                                                                      & =  \lim_{r \to \infty} \inf\frac{1}{r^2}\, C, \, \operatorname{for}\, \operatorname{some}\, \operatorname{constant}\, C > 0 \\
& = 0 < \infty. \end{aligned} \] 
Definition \ref{def2.1} and the first assertion of Theorems \ref{T:1.1} complete the proof.\hskip.4in  $\square$

\vskip.05in
\noindent{\it Proof of Theorem \ref{T:1.3}}. Let $\phi: N=N_1\times_{f_2} N_2 \times \cdots \times _{f_k} N_k \to M$ be an isometric immersion of a multiply warped product $N = N_1\times_{f_2} N_2 \times \cdots \times _{f_k} N_k$ in a Riemannian manifold $M$. If the mean curvature $H$ of $N$ in $M$ satisfies \eqref{1.7} on $N$, then it follows from the inequality \eqref{1.4} of Theorem D and \eqref{1.7} that
$ \sum _{j=2}^{k} n_j  \frac{\Delta f_j}{f_j}  \ge 0.$  Hence, after applying the Average Method in PDE stated in Proposition \ref{P:2.1}(1),  we  conclude  that  some $f_i,\, 2\leq i\leq k,$ could be constant {\it or} $f_i$ would be $2$-imbalanced.\hfill $\square$

\vskip.05in
\noindent{\it Proof of Corollary \ref{C:1.1}}.  Suppose contrary, such an immersion would violate \eqref{1.5} and hence contradicts Theorem \ref{T:1.1}. \hfill $\square$

\vskip.05in
\noindent{\it Proof of Corollary \ref{C:1.2}}. Statement $(1)$ of Corollary \ref{C:1.2} follows from Theorem \ref{T:1.3} and Proposition \ref{P:2.1}(2). 

 In view of \eqref{1.7} and Proposition \ref{P:2.1}(2), we have 
\begin{equation}\notag 0 \le- H^2  - \frac{2kn_1(n - n_1) c }{n^2(k - 1)} \le \sum _{j = 2}^ k  n_j\frac{\Delta f_j}{f_j} = 0.\end{equation} 
Therefore we obtain $H = c = 0$, which implies statement (2). 

Statement (3) follows immediately from Corollary \ref{C:1.5}. \hfill $\square$

\vskip.05in
\noindent{\it Proof of Corollary \ref{C:1.3}}.  To prove statement (1),  let us suppose contrary. Then it follows from Theorem D that 
\begin{equation} \label{3.3} 0 <- H^2  - \frac{2kn_1(n - n_1) \max \tilde K }{n^2(k - 1)} \le \sum _{j = 2}^ k  n_j \frac{\Delta f_j}{f_j}.\end{equation} 
Now, by Theorem \ref{T:2.1}(1), \eqref{3.3}  implies the constancy of   $f_i$ for some $2 \le i \le k\, .$ Thus 
\begin{equation} 0 < \sum _{j\ne i}  \frac{\Delta f_j}{f_j} \label{3.4}. \end{equation}
Therefore, after applying Proposition \ref{P:2.1}(2) to \eqref{3.4} we obtain the constancy of $f_2, \cdots , f_k$, which leads to  
$0 < 0\, ,$ a contradiction.

For statement (2), let us suppose contrary. Then inequality \eqref{3.2} would be true. Hence by Proposition \ref{P:2.1}(1), we would have the constancy of   $f_i$ for some $2 \le i \le k\, .$ Thus 
\begin{equation} 0 \le \sum _{j\ne i} n_j  \frac{\Delta f_j}{f_j} \label{3.5}. \end{equation}
Now, applying Proposition \ref{P:2.1}(2) shows the constancy of $f_2, \cdots , f_k$. 
So, it follows from Theorem D that
 $\phi$ is mixed totally geodesic and hence, by Moore's lemma \cite {M},  we conclude that $$\phi = (\phi _1, \cdots , \phi _k) : N = N_1\times_{f_2} N_2 \times \cdots \times _{f_k} N_k \to {\mathbb E}^{m_1} \times \cdots \times {\mathbb E}^{m_k}= {\mathbb E}^m$$  
 is a product minimal immersion, which contradicts to the fact that there is no compact minimal
submanifold $N_k$ in the Euclidean space ${\mathbb E}^{m_k}$. \hfill $\square$

\vskip.05in
\noindent{\it Proof of Corollary \ref{C:1.4}}.  Follows at once from Corollary \ref{C:1.3} and the fact that every $L^{q}$ function with $q  > 1$ on $N_1$ is $2$-balanced for the same $q > 1$ on $N_1.$ \hfill $\square$

\vskip.05in
\noindent{\it Proof of Corollary \ref{C:1.5}}.  In view of Theorem D and $H = c = 0\, ,$ we have 
 \[ 0 = - H^2  - \frac{2kn_1(n - n_1) c }{n^2(k - 1)} = \sum _{j = 2}^ k  n_j \frac{\Delta f_j}{f_j} = 0.\] 
Now assertion follows from Theorem E \cite[N\"olker's Theorem]{N}. \hfill $\square$

\begin{remark} \label{R:3.1}{\rm The following two examples show that Theorem \ref{T:1.1} is false if either $f_j$ is constant or $f_j$ is 2-imbalance for every $q_j> 1$.} 
\vskip.1in

\noindent{\bf Example 3.1.} {\rm Let $N_1,\ldots,N_k$ be $k$ copies of the real line ${\bf R}$ and let us put  $f_j=1$ for $j=2,\ldots,k$. Then  $N=N_1 \times_{1} N_2 \times\cdots \times_{1} N_k$ is the Euclidean $k$-space $\mathbb E^k$. Clearly, for a totally geodesic immersion of $N$ into $\mathbb E^{k+1}$, inequality \eqref{1.5} is false.}
\vskip.1in

\noindent{\bf Example 3.2.}  {\rm Let $N_1=\{x\in {\bf R}: x>0\}$ and $N_2,\ldots,N_k$ be $k-1$ copies of  ${\bf R}$. If we put  $f_2=\cdots=f_k=x$, then each $f_j$ is 2-imbalance for every $q_j>1$
and $N=N_1 \times_{x} N_2 \times\cdots \times_{x} N_k$ is an open subset of  $\mathbb E^k$. Again, for a totally geodesic immersion of $N$ into $\mathbb E^{k+1}$, inequality \eqref{1.5} is false.}
\end{remark}

\begin{remark} \label{R:3.2} {\rm Theorems \ref{T:1.1} and \ref{T:1.2} are sharp in the sense that inequality \eqref{1.5} and \eqref{1.6} are false if either $f_j$ were constant or $f_j$ were 2-imbalanced for every $q_j>1$ (For details, we refer to Remark \ref{R:3.1}, and Examples 3.1--3.2 above).  Furthermore, Theorem \ref{T:1.3} shows that inequality \eqref{1.5} (resp., \eqref{1.6}) is best possible for Theorem \ref{T:1.1} (resp., for Theorem \ref{T:1.2}).}
\end{remark}

\section{\bf Multiply Warped Product  Manifolds into complex or quaternionic space forms}

A submanifold $N$ of a K\"ahler manifold $M$ is said to be {\it totally real} if the almost complex structure $J$ of $M$ carries each tangent space of $N$ into its corresponding normal space (cf. \cite{CO,c1}). Similarly, one has the notion of totally real submanifolds in quaternionic K\"ahler manifolds (cf. \cite{CH}).

B.-Y. Chen and F. Dillen proved
\vskip.1in

{\bf Theorem F} \cite{CD} {\it Let $\phi : N_1\times_{f_2} N_2 \times \cdots \times _{f_k} N_k  \to \tilde M^m (4c)$ be a totally real isometric 
immersion  of the multiply warped product $N = N_1\times_{f_2} N_2 \times \cdots \times _{f_k} N_k $ into a  complex space form of constant holomorphic sectional curvature $4c$ or in a quaternionic space form of constant  quaternionic sectional curvature $4c$. Then 
\begin{equation}\label{4.1} -\frac{n^2}{4} H^2-n_1(n-n_1)c \leq \sum_{j=2}^k n_j \frac{\Delta f_j}{f_j} ,\quad n=\sum_{i=1}^k n_i.\end{equation}}
\vskip.1in

By applying the same techniques, i.e., Theorem \ref{T:2.1} (Warping Functions Growth Estimates)  and Proposition \ref{P:2.1} (An Average Method in PDE), we also have the following results.

\begin{theorem} \label{T:4.1}
If for each $j, 2 \le j \le k$, $f_j$ is nonconstant and \emph{$2$-balanced} for some $q_j > 1$,
then for any multiply warped product $N = N_1\times_{f_2} N_2 \times \cdots \times _{f_k} N_k$ totally real isometrically immersed in a complex space form of constant holomorphic sectional curvature $4c$ or in a quaternionic space form of constant  quaternionic sectional curvature $4c$, the mean
curvature $H$ of $\phi$ satisfies
\begin{equation}\label{4.2}H^2 > \frac{- 4n_1(n - n_1)c}{n^2}\end{equation} at some points. 

In particular, if each $f_j $ is nonconstant and in $L^{q_j}$ for some $q_j > 1$, then \eqref{4.2} holds at some points.\end{theorem}

As applications of Theorem \ref{T:4.1} , we have the following.

\begin{corollary} \label{C:4.1} If for each $j, 2 \le j \le k$, $f_j$ is nonconstant and \emph{$2$-balanced} for some $q_j > 1$, then there does not exist a totally real minimal immersion of any multiply warped product $N = N_1\times_{f_2} N_2 \times \cdots \times _{f_k} N_k$ into a complex space form of constant holomorphic sectional curvature $4c\leq 0$ or into a quaternionic space form of constant  quaternionic sectional curvature $4c\leq 0$. 

In particular,  if each $f_j $ is nonconstant and in $L^{q_j}$ for some $q_j > 1$, then there does not exist an isometric minimal immersion of $N = N_1\times_{f_2} N_2 \times \cdots \times _{f_k} N_k$ into $\tilde M^m (0)$.
\end{corollary}

Another application of Theorem \ref{T:4.1} is the following dichotomy.  

\begin{corollary} \label{C:4.2}
Suppose the squared mean curvature of an isometric immersion of a multiply warped product $N = N_1\times_{f_2} N_2 \times \cdots \times _{f_k} N_k$ into a Riemannian manifold $R^m(c)$ of constant sectional curvature $c$ satisfies
\begin{align}\label{4.3}H^2 \leq \frac{- 4n_1(n - n_1)c}{n^2}\end{align} everywhere on $N$. Then there exists an integer $i, 2 \le i \le k$, such that the warping function $f_i$ is either a constant or for
every $q_i > 1\, ,$ $f_i$ has \emph{$2$-imbalanced} growth. 
\end{corollary}

By applying the Growth Estimates in Theorem \ref{T:2.1} and an Average Method in PDE in Proposition \ref{P:2.1}, we have the following.

\begin{corollary} \label{C:4.3} Suppose the squared mean curvature of a totally real isometric immersion $\phi$ of a multiply warped product $N = N_1\times_{f_2} N_2 \times \cdots \times _{f_k} N_k$ into  a complex space form of constant holomorphic sectional curvature $4c$ or a quaternionic space form of constant  quaternionic sectional curvature $4c$ satisfies \eqref{4.3} everywhere on $N$.  If for each $j, 2 \le j \le k$, $f_j$ is \emph{$2$-balanced} for some $q_j > 1$,
then  
\begin{enumerate}
\item[{(1)}]  Every  warping function $f_j, 2 \le j \le k,$ is constant. 

\item[{(2)}] The isometric immersion $\phi$ is a minimal immersion into $\tilde M^m (0)$.  
\end{enumerate}\end{corollary}

\begin{corollary} \label{C:4.4} Let each $f_j\, ,$ $2 \le j \le k,$  be  \emph{$2$-balanced} for some $q_j > 1\, .$ Then
\begin{enumerate}
\item[{(1)}]  Every multiply warped product $N = N_1\times_{f_2} N_2 \times \cdots \times _{f_k} N_k$ does not admit  an isometrically totally real minimal immersion into any complex space form of negative constant holomorphic sectional curvature $4c$ or a quaternionic space form of negative constant  quaternionic sectional curvature $4c$.

\item[{(2)}]  If $N_1$ is compact,  then $N = N_1\times_{f_2} N_2 \times \cdots \times _{f_k} N_k$ does not admit  an isometrically totally real minimal immersion into $\tilde M^m (0)$. 
\end{enumerate}\end{corollary}

We state a special case of Corollary \ref{C:4.4} as follows.

\begin{corollary} \label{C:4.5}   If each $f_j\, ,$ $2 \le j \le k,$ is  in $L^{q_j}$ for some $q_j > 1$, then we have:
\begin{enumerate}
\item[{(1)}] Every multiply warped product $N = N_1\times_{f_2} N_2 \times \cdots \times _{f_k} N_k$ does not 	admit  an isometrically totally real minimal immersion into a complex space form of negative constant holomorphic sectional curvature $4c$ or into a quaternionic space form of negative constant  quaternionic sectional curvature $4c$. 
\item[{(2)}] If $N_1$ is compact,  then  $N_1\times_{f_2} N_2 \times \cdots \times _{f_k} N_k$ does not admit  an isometrically totally real minimal immersion into $\tilde M^m (0)$. 
\end{enumerate}
\end{corollary}

Since the proofs of these results can be done in the same way as in section 3, we omit their proofs.

\section{\bf  Doubly Warped Products}

Doubly warped products are natural generalization of (ordinary) warped products. 

\begin{definition} A {\it doubly warped product} of  Riemannian manifolds $(N_1,g_1)$ and $(N_2,g_2)$ is a product manifold  $ {}_{f_2}\hskip-.0in N_1\times_{f_1}\hskip-.04in  N_2$ equipped with metric $g = f_2^2 g_1 \oplus f_1^2 g_2\, ,$ where $f_1 : N_1 \to \mathbb R^{+}$ and $f_2: N_2 \to \mathbb R^{+}$ are positive-valued smooth functions.
\end{definition}

As an extension of Theorem A from \cite{CW}, A. Olteanu proved the following.
\vskip.1in

\noindent
{\bf Theorem G} \cite {O}  {\it Let $\phi : {}_{f_2}\hskip-.0in N_1\times_{f_1}\hskip-.04in  N_2 \to M $ be an isometric 
immersion  of a doubly warped product $ {}_{f_2}\hskip-.0in N_1\times_{f_1}\hskip-.04in  N_2$ into an arbitrary Riemannian manifold $M$.
We have}
\begin{equation}\label{5.1} - \frac{(n_1+n_2)^2}{4} H^2- n_1 n_2\max \tilde{K} \leq n_2 \frac{\Delta_1 f_1}{f_1}+n_1 \frac{\Delta_2 f_2}{f_2},\end{equation}
{\it where $n_i = \dim N_i$ and $\Delta_i$ is the Laplacian of $N_i\, ,$ for $i= 1, 2.$ 

The equality sign holds identically if and only if the following conditions hold: }
\begin{enumerate}
\item[{\rm i)}]  {\it $\phi$ is mixed totally geodesic such that
$\operatorname{trace}\, h _1 =  \operatorname{trace}\, h _2\, $; }

\noindent
\item[{\rm ii)}]  {\it At each point $x=(x_1,x_2) \in N\, , $ $\tilde K$ satisfies $\tilde K(u, v) = \max K(x)$ for each unit vector $u \in T^1 _{x_1} N_1$ and every  $v \in T^1 _{x_2} N_2.$ }
\end{enumerate}
\vskip.1in

Similarly, by applying the same techniques via Theorem \ref{T:2.1} (Warping Functions Growth Estimates)  and Proposition \ref{P:2.1} (An Average Method in PDE), we also have the following.

\begin{theorem} \label{T:5.1} If $f_1,f_2$ are nonconstant and \emph{$2$-balanced} for some $q_1,q_2 > 1$,
then for any  isometric  immersion of a doubly warped product $\phi : {}_{f_2}\hskip-.0in N_1\times_{f_1}\hskip-.04in  N_2$  into a
Riemannian manifold $M$, the mean
curvature $H$ of $\phi$ satisfies
\begin{equation}\label{5.2}H^2 > \frac{- 4n_1 n_2}{(n_1 + n_2)^2}\max \tilde K\end{equation} at some points. 

In particular, if each $f_j $ is nonconstant and in $L^{q_j}$ for some $q_j > 1$, then \eqref{5.2} holds at some points.\end{theorem}

The proof of this theorem is similar to the proof of Theorem \ref{T:2.1}. However, because doubly warped products are somewhat different from ordinary warped products, we provide the proof of Theorem 5.1 as follows. 

\vskip.1in
\noindent{\it Proof of Theorem \ref{T:5.1}.} Suppose contrary to \eqref{5.2}, i.e., there were an isometric immersion $\phi$ whose mean curvature $H$ satisfying
\begin{equation}\notag H^2 \leq \frac{- 4n_1 n_2}{(n_1 + n_2)^2}\max \tilde K\end{equation}  everywhere on $N$, which gives
\begin{equation}\label{5.3}0 \le - \frac {(n_1 + n_2)^2}{4} H^2 - n_1 n_2 \max \tilde{K}. \end{equation}

On the other hand, Theorem G would imply
\begin{equation}\notag - \frac{(n_1+n_2)^2}{4} H^2- n_1 n_2\max \tilde{K} \leq n_2 \frac{\Delta_1 f_1}{f_1}+n_1 \frac{\Delta_2 f_2}{f_2}.\end{equation}
After combining this with \eqref{5.3}, we find
\begin{equation}
\begin{aligned} n_2 \frac{\Delta_1 f_1}{f_1}+n_1 \frac{\Delta_2 f_2}{f_2} & \ge - \frac{(n_1+n_2)^2}{4} H^2- n_1 n_2\max \tilde{K} \ge 0.
 \end{aligned} \label{5.4} 
 \end{equation}
After applying the Average Method in PDE stated in Proposition \ref{P:2.1}(1),   \eqref{5.4} shows that 
 some $f_i,\,i=1,2,$ could be constant {\it or} $f_i$ would be $2$-imbalanced for every $q_i > 1$. This contradicts the assumption that
$f_i$ is nonconstant {\it and}  $2$-balanced for some $q_i > 1.$ 

The last assertion follows from Definition \ref{def2.1}, the first assertion of Theorems \ref{T:5.1} and the fact that every $L^q$ function has \emph{$2$-{finite}}, \emph{$2$-mild}, \emph{$2$-obtuse}, \emph{$2$-moderate}, \emph{$2$-small} growth for the same $q$ (cf. \cite [Proposition 2.3] {WLW}). 
\hfill $\square$

In particular, if the ambient space $M$ is of constant sectional curvature $c\leq 0$,  then Theorem \ref{T:5.1} reduces to the following.

\begin{theorem} \label{T:5.2} If $f_1,f_2$ are nonconstant and \emph{$2$-balanced} for some $q_1,q_2 > 1$,
then for any  isometric  immersion of a doubly warped product $\phi : {}_{f_2}\hskip-.0in N_1\times_{f_1}\hskip-.04in  N_2$  into a
Riemannian $m$-manifold $R^m(c)$ of constant curvature $c\leq 0$, the mean
curvature $H$ of $\phi$ satisfies
\begin{equation}\label{5.5}H^2 > \frac{-4n_1 n_2}{(n_1 + n_2)^2}c\end{equation} at some points. 

In particular, if each $f_j $ is nonconstant and in $L^{q_j}$ for some $q_j > 1$, then \eqref{5.5} holds at some points.
\end{theorem}

Also, the following are easy consequences of Theorem \ref{T:5.1}.

\begin{corollary} \label{C:5.1} If for each $j\, (j=1,2)$, $f_j$ is nonconstant and \emph{$2$-balanced} for some $q_j$,
then there does not exist an  isometric minimal 
immersion of any doubly warped product $\phi :  {}_{f_2}\hskip-.0in N_1\times_{f_1}\hskip-.04in  N_2$  into any
 negatively curved Riemannian manifold.  
 
  In particular,  if each $f_j $ is nonconstant and in $L^{q_j}$ for some $q_j > 1$, then there does not exist isometric minimal immersion of  $N = N_1\times_{f_2} N_2 \times \cdots \times _{f_k} N_k$ into a Euclidean space.
\end{corollary}

As another easy applications of Theorem \ref{T:5.1}, we have the following dichotomy.  

\begin{corollary} \label{C:5.2} Suppose the squared mean curvature of the isometric immersion of a doubly warped product $\phi : _{f_2}N_1\times_{f_1} N_2 $  into a Riemannian manifold satisfies
\begin{equation} H^2 \leq \frac{- 4n_1 n_2}{(n_1 + n_2)^2}\max \tilde K\end{equation} everywhere on $N$. Then there exists an integer $i, 1 \le i \le 2,$ such that the warping function $f_i$ is either a constant or for
every $q_i > 1\, ,$ $f_i$ has \emph{$2$-imbalanced} growth. 
\end{corollary}

Analogously, by applying the growth estimates and the average method in PDE as before, we have the following Liouville property and a characterization result.

\begin{corollary} \label{C:5.3} Suppose the squared mean curvature of an
isometric immersion $\phi$ of a doubly warped product ${}_{f_2}\hskip-.0in N_1\times_{f_1}\hskip-.04in  N_2$  into a Riemannian manifold satisfies \eqref{5.3} on $N$.  If for each $j\, (j=1,2)$, $f_j$ is \emph{$2$-balanced} for some $q_j > 1$, then we have:
\begin{enumerate}
\item[{\rm (1)}]  Every  warping function $f_j, 1 \le j \le 2$, is constant. 
\item[{\rm (2)}]  The isometric immersion $\phi$ is a minimal immersion into a Euclidean space.  
\end{enumerate}\end{corollary}

\begin{corollary} \label{C:5.4} Let each $f_j\, ,$ $1 \le j \le 2$  be  \emph{$2$-balanced} for some $q_j > 1\, .$ Then we have:
\begin{enumerate}
\item[{\rm (1)}] Every doubly warped product $ {}_{f_2}\hskip-.0in N_1\times_{f_1}\hskip-.04in  N_2$ does not admit  an isometric minimal immersion into any negatively curved Riemannian manifold.

\item[{\rm (2)}] If $N_2$ is compact,  then $ _{f_2}N_1\times_{f_1} N_2 $ does not admit  an isometric minimal immersion into a Euclidean space. 
\end{enumerate}\end{corollary}

We state a special case of Corollary \ref{C:5.5}:

\begin{corollary} \label{C:5.5}   If each $f_j\,(j=1,2)$ is  in $L^{q_j}$ for some $q_j > 1$, then we have:
\begin{enumerate}
\item[{\rm (1)}]  Every doubly warped product ${}_{f_2}\hskip-.0in N_1\times_{f_1}\hskip-.04in  N_2$ does not admit  an isometric minimal immersion into any negatively curved Riemannian manifold. 
\item[{\rm (2)}]  If $N_2$ is compact,  then ${}_{f_2}\hskip-.0in N_1\times_{f_1}\hskip-.04in  N_2$ does not admit  an isometric minimal immersion into a Euclidean space. 
\end{enumerate}\end{corollary}

Since Corollaries \ref{C:5.1}-\ref{C:5.5} can be in the same way as the proofs of Corollary \ref{C:1.1}, Theorem \ref{T:1.3}, Corollary \ref{C:1.2} (1)\&(2), Corollary \ref{C:1.3} and Corollary \ref{C:1.4}, we omit their proofs.

\vfill\eject

\vskip.2in
\noindent {B.-Y. Chen, Department of Mathematics,
    Michigan State University, East Lansing, Michigan 48824--1027, U.S.A.}
\\
{\it E-mail address}: {\tt chenb@msu.edu}
\vskip.2in

\noindent {S. W. Wei, Department of Mathematics, University of Oklahoma, Norman, Oklahoma 73019-0315, U.S.A.}\\
{\it E-mail address}: {\tt wwei@ou.edu}

\end{document}